\documentclass[preprint,12pt]{elsarticle}

\usepackage{amssymb}
\usepackage{amsfonts}
\usepackage{amsmath}
\usepackage{latexsym,epsfig}
\usepackage{graphicx} 
\usepackage[usenames,dvipsnames]{color} 
\usepackage{array} 
\usepackage{color} 
\usepackage{verbatim}
\usepackage[ruled,linesnumbered]{algorithm2e}  
\usepackage{url} 
\usepackage{enumerate} 
\usepackage{rotating}
\usepackage{hyperref} 
\usepackage{afterpage,lscape} 
\usepackage[misc,geometry]{ifsym} 
\usepackage[shortlabels]{enumitem} 
\usepackage{wasysym} 
\setlist[enumerate]{after={\bigskip}}
\setlist[itemize]{after={\bigskip}} 
\allowdisplaybreaks 
 
\newcommand{\EndProof}{\hspace{\stretch{1}} $\Box$}

\newcommand{\Z}{\mathbb{Z}}

\newcommand{\N}{\mathbb{N}}
\newcommand{\M}{\mathbb{M}}
\newcommand{\Semi}{\mathbb{S}}

\newcommand{\MREP}{\textsc{MinRep}} 
\newcommand{\MC}{\textsc{MinCost}} 
\newcommand{\GREP}{\textsc{GreedyRep}} 
\newcommand{\GRC}{\textsc{GreedyCost}}
\newcommand{\QGR}{\textsc{GenGreedyRep}} 
\newcommand{\QGC}{\textsc{GenGreedyCost}} 

\newtheorem{definition}{Definition}
\newtheorem{example}{Example}
\newtheorem{lemma}{Lemma}
\newtheorem{proposition}{Proposition}
\newtheorem{theorem}{Theorem}
\newtheorem{observation}{Observation}

\journal{Communications in Algebra, accepted}

\begin{document}


\begin{frontmatter}
	
	\title{Greedy Sets and Greedy Numerical Semigroups}
	
	\author{ Hebert P\'erez-Ros\'es\footnote{\Letter  Corresponding author: hebert.perez@urv.cat} \\ Jos\'e Miguel Serradilla-Merinero \\ Maria Bras-Amor\'os }  
	\affiliation{organization={Dept. of Computer Science and Mathematics \\ Universitat Rovira i Virgili}, 
		addressline={Avda. Paisos Catalans 26}, 
		city={Tarragona}, 
		postcode={43007}, 
		state={Catalonia}, 
		country={Spain}} 
	
	
	\begin{abstract}
		Motivated by the change-making problem, we extend the notion of greediness to sets of positive integers not containing the element $1$, and from there to numerical semigroups. We provide an algorithm to determine if a given set (not necessarily containing the number $1$) is greedy. We also give specific conditions for sets of cardinality three, and we prove that numerical semigroups generated by three consecutive integers are greedy.   
	\end{abstract}

	
\end{frontmatter}

\section{Greedy sets and greedy numerical semigroups}
\label{sec:intro} 

In the money-changing problem, or change-making problem, we have a set of integer coin denominations $S = \{ s_1, s_2, \ldots, s_t \}$, with $0 < s_1 < \ldots < s_t$. We also have a target amount $k$, and the goal is to make $k$ using as few coins as possible. In order to be able to represent any amount $k$ we usually impose the additional requirement that $s_1 = 1$, although in this paper we will drop that requirement later. 

More formally, we are looking for a \emph{payment vector} $\textbf{a} = (a_1, \ldots, a_t)$, such that 

\begin{align} 
	\label{eq:cond1} & a_i \in \N_0, \mbox{ for all } i = 1, \ldots, t \\
	\label{eq:cond2} & \sum_{i=1}^{t} a_i s_i = k, \\
	\label{eq:cond3} & \sum_{i=1}^{t} a_i \mbox{ is minimal,}  
\end{align} 
where $\N_0$ denotes the set of nonnegative integers. 

\begin{definition}
	\label{def:mincost} 
	Given a set of denominations $S=\{ s_1, s_2, \ldots, s_t \}$ and a given representable amount $k > 0$, a \emph{minimal payment vector}, or \emph{minimal representation} of $k$ with respect to $S$, denoted $\mbox{\MREP}_S(k)$ is a payment vector $\textbf{a} = (a_1, a_2, \ldots, a_t)$ that satisfies Conditions (\ref{eq:cond1}), (\ref{eq:cond2}) and (\ref{eq:cond3}). If $\textbf{a}$ is a minimal representation of $k$, then   $\displaystyle \mbox{\MC}_S(k) = \sum_{i=1}^{t} a_i$. 
\end{definition} 

\noindent Note that $\mbox{\MREP}_S(k)$ is not necessarily unique, but $\mbox{\MC}_S(k)$ is. Another simple property of the function $\mbox{\MC}_S(k)$ is given by

\begin{proposition}[Triangle inequality] Let $h$ and $k$ be two representable numbers. Then 
    \label{prop:triangularinequality1} 
    \begin{equation}
        \label{eq:triangularinequality1}
        \mbox{\MC}_S(h+k) \leq \mbox{\MC}_S(h) + \mbox{\MC}_S(k) 
    \end{equation}
\end{proposition}

\noindent \textbf{Proof:} Let $\textbf{h} = \mbox{\MREP}_S(h) = (h_1, h_2, \ldots, h_t)$ and $\textbf{k} = \mbox{\MREP}_S(k) = (k_1, k_2, \ldots, k_t)$ be the optimal representations of $h$ and $k$, respectively. We have 
\begin{displaymath}
    \sum_{i=1}^t h_i s_i = h, \qquad \sum_{i=1}^t k_i s_i = k,
\end{displaymath}
Now we just have to check that $\textbf{a} = (h_1 + k_1, h_2 + k_2, \ldots, h_t + k_t)$ is a representation of $h+k$ (not necessarily  optimal) since 
\begin{displaymath}
    \sum_{i=1}^t (h_i + k_i) s_i = h+k. 
\end{displaymath}

\EndProof

\vspace{2mm}
A traditional approach for solving the change-making problem is the \emph{greedy approach}, which proceeds by first choosing the coin of the largest possible denomination, subtracting it from the target amount, and then applying the same algorithm to the remainder. The  method is formally described in Algorithm \ref{alg:greedypayment} (see \cite{PR21}). 
\\\\
	\begin{algorithm}[H]
		\SetKwInOut{Input}{Input}
		\SetKwInOut{Output}{Output}
		\SetKw{DownTo}{downto}
		\vspace{.2cm}
		\Input{The set of denominations $S = \{ 1, s_2, \ldots, s_t \}$, with $1 < s_2 < \ldots < s_t$, and a target amount $k > 0$.} 
		\Output{Payment vector $\textbf{a} = (a_1, a_2, \ldots, a_t)$.} 
		\vspace{.2cm}
		
		\For{ $i$:= $t$ \DownTo $1$} 
		{
			$a_i$ := $k$ div $s_i$\;
			$k$ :=  $k$ mod $s_i$\;
            \If{ $k=0$ }{ \Return{ $\textbf{a}$ } } 
		} 
		\caption{GREEDY PAYMENT METHOD}
		\label{alg:greedypayment}  
	\end{algorithm}
\vspace{3mm} 

Now let's suppose that the set of coin denominations $S = \{ s_1, s_2, \ldots, s_t \}$ is such that $1 < s_1 < \ldots < s_t$. This obviously entails some difficulties, the first one being that not all integers are representable (whatever the algorithm we choose for finding the representation). 

\begin{definition}
    \label{def:greedy-representable} 
    Given a set $S = \{ s_1, s_2, \ldots, s_t \}$ of denominations, such that $1 < s_1 < \ldots < s_t$, and a target amount $k$, we say that $k$ is \emph{representable} if there exists a payment vector $\textbf{a} = (a_1, a_2, \ldots, a_t)$ satisfying Conditions (\ref{eq:cond1}) and (\ref{eq:cond2}). We say that $k$ is \emph{greedy-representable} if there exists a payment vector $\textbf{a} = (a_1, a_2, \ldots, a_t)$ satisfying Conditions (\ref{eq:cond1}) and (\ref{eq:cond2}),  that is obtained by Algorithm \ref{alg:greedypayment}. 

    \noindent Additionally, if $k$ is greedy-representable, then the \emph{greedy payment vector}, or \emph{greedy representation} of $k$ with respect to $S$, is the payment vector $\textbf{a} = \mbox{\GREP}_S(k) = (a_1, a_2, \ldots, a_t)$ produced by Algorithm \ref{alg:greedypayment}, and 
    \begin{displaymath}
        \mbox{\GRC}_S(k) = \sum_{i=1}^{t} a_i.
    \end{displaymath} 
\end{definition} 


Now, if $\gcd(s_1, s_2, \ldots, s_t) = 1$ then $S$ generates a numerical semigroup $\Semi$, and there exists an integer $F(\Semi)$ (the Frobenius number of $\Semi$), such that any integer $k > F(\Semi)$ is representable. In this context the denominations $s_1, s_2, \ldots, s_t$ are called the \emph{generators} of $\Semi$, and we write $\Semi = \langle S \rangle$. We also denote by $G(\Semi)$ the set of gaps of $\Semi$. For the basic concepts and results of numerical semigroups see \cite{RoGar09}.  

So, from now on we may assume that $\gcd(s_1, s_2, \ldots, s_t) = 1$ and $k > F(\Semi)$, so we can guarantee that $k$ is representable. Note that in this new scenario the definition of $\mbox{\MREP}_S(k)$ and $\mbox{\MC}_S(k)$ remains unchanged. The following lemma is adapted from \cite{KoZaks94}, and will be needed later: 
\begin{lemma}
\label{lemma:kozaks}
    Let $\Semi = \langle S \rangle$ be a numerical semigroup generated by the set $S = \{ s_1, \ldots, s_t \}$. Additionally, let $k \in \Semi$ and $s_i$ such that $k-s_i \in \Semi$. Then 
    \begin{equation}
        \label{eq:lemmakozen}
        \mbox{\MC}_S(k) \leq \mbox{\MC}_S(k-s_i)+1, 
    \end{equation}
    with equality holding if, and only if, there exists an optimal representation of $k$ that uses the generator $s_i$. 
\end{lemma}
\vspace{4mm} 
\noindent \textbf{Proof:} It is clear that Eq. \ref{eq:lemmakozen} holds. If additionally $\mbox{\MC}_S(k) = \mbox{\MC}_S(k-s_i)+1$, then the representation of $k$ so obtained is optimal and uses the generator $s_i$. Conversely, given an optimal representation of $k$ that uses $s_i$, we can obtain a representation of $k-s_i$ of cost $\mbox{\MC}_S(k)-1$ by subtracting one from the coefficient of $s_i$ in $\mbox{\MREP}_S(k)$, and thus Eq. \ref{eq:lemmakozen} implies that this representation is optimal. 

\EndProof 

\vspace{2mm}
Regardless of the fact that our $k$ is representable in the new scenario (i.e. with $s_1 > 1$), we cannot still guarantee that $k$ be greedy-representable in the sense of Definition \ref{def:greedy-representable}. In fact, we have the following statement:  

\begin{proposition}
    \label{prop:badness1}
    Let $\Semi = \langle S \rangle$ be a numerical semigroup generated by $S = \{ s_1, s_2, \ldots, s_t \}$. Then there exists an infinite set of integers $\M \subseteq \Semi$ that are representable but not greedy-representable. 
\end{proposition} 

\noindent \textbf{Proof:} Let $g \in G$ be any gap of $\Semi$, and let $\M_g = \{ g + j s_t : j \in \N \}$. It is easy to check that $\M_g \subseteq \M$.   

\EndProof 
\\\\ 
Let us denote by $\mbox{Ap}(\Semi,n)$ the Ap\'ery set of $\Semi$ with respect to $n$, i.e. 
\begin{displaymath}
    \mbox{Ap}(\Semi,n) = \{ s \in \Semi : s-n \notin \Semi \}. 
\end{displaymath}
Then $\mbox{Ap}(\Semi,s_t) \subseteq \M$, where $\M$ is the set mentioned in Proposition \ref{prop:badness1}.  
\\\\ 
This fact leads us to consider another strategy for representing numbers, which we call \emph{generalized greedy}, and is formalized in Algorithm \ref{alg:quasigreedypayment}. 
\\\\ 
	\begin{algorithm}[H] 
		\SetKwInOut{Input}{Input}
		\SetKwInOut{Output}{Output}
		\SetKw{DownTo}{downto}
		\vspace{.2cm}
		\Input{The set of denominations $S = \{ s_1, s_2, \ldots, s_t \}$, with $1 \leq s_1 < s_2 < \ldots < s_t$, $\gcd(s_1, s_2, \ldots, s_t) = 1$, and an element $k \in \langle S \rangle$, $k > 0$.} 
		\Output{Generalised greedy representation vector $\textbf{a} = (a_1, a_2, \ldots, a_t)$.} 
		\vspace{.2cm}
		
		\For{ $i$:= $t$ \DownTo $1$} 
		{ 
            \While{ $k-s_i \in \langle S \rangle$ }
            { 
			 $a_i$ := $a_i + 1$\;
			 $k$ :=  $k - s_i$\;  
              \If{ $k=0$ }{ \Return{ $\textbf{a}$ }\; } 
            }
		} 
		\caption{GENERALISED GREEDY REPRESENTATION}
		\label{alg:quasigreedypayment}  
	\end{algorithm}
\vspace{3mm} 
Note that the generalized greedy algorithm uses the greedy strategy as much as possible. It chooses the largest possible coin denomination until it becomes apparent that this greedy strategy will lead to a dead end, and then the algorithm proceeds to the next largest denomination. 

Note also that Algorithm \ref{alg:quasigreedypayment} always terminates. Indeed, the algorithm always moves from one representable element $k \in \langle S \rangle$ to a smaller representable element $k'$, eventually arriving at zero. The work of Algorithm \ref{alg:quasigreedypayment} is illustrated in Example \ref{ex:algorithmic-subtract}. 

\begin{example}
    \label{ex:algorithmic-subtract}
    Let $S = \{ 3, 7, 10 \}$ be a set of denominations and $\Semi = \langle S \rangle$. If we apply Algorithm \ref{alg:quasigreedypayment} to $k=28$ we observe the steps summarized in Table \ref{tab:execution-algorithm-subtract}. We have ommitted some steps for the sake of economy, but the reader should have no difficulty in completing them. The algorithm terminates and outputs the representation vector $(6, 0, 1)$. Notice that $\mbox{\MREP}_S(28) = (0, 4, 0)$. 
    \begin{table}[h]
    \footnotesize 
    \centering
        \begin{tabular}{|c|c|c|c|c|c|c|}
            \hline
            Index $i$ & Target $k$ & Gen. $s_i$ & $a_i$ & $k - s_i$ & Test & Actions \\ \hline
            3  & 28  & 10 & $a_3 = 0$ & 18 & $18 \in \Semi$ & $k \leftarrow 18$, $a_3 \leftarrow a_3 + 1$  \\ 
            3  & 18  & 10 & $a_3 = 1$ & 8 & $8 \notin \Semi$ & $i \leftarrow i - 1$  \\  
            2  & 18  & 7 & $a_2 = 0$ & 11 & $11 \notin \Semi$ & $i \leftarrow i - 1$  \\ 
            1  & 18  & 3 & $a_1 = 0$ & 15 & $15 \in \Semi$ & $k \leftarrow 15$, $a_1 \leftarrow a_1 + 1$  \\ 
            1  & 15  & 3 & $a_1 = 1$ & 12 & $12 \in \Semi$ & $k \leftarrow 12$, $a_1 \leftarrow a_1 + 1$  \\
            \vdots & \vdots & \vdots & \vdots & \vdots & \vdots & \vdots \\ 
            1  & 3  & 3  & $a_1 = 5$  & 0  & $0 \in \Semi$ & $k \leftarrow 0$, $a_1 \leftarrow a_1 + 1$ \\ 
             & & & & & & Return $(6, 0, 1)$ \\ 
            \hline
        \end{tabular} 
    \caption{Execution of Algorithm \ref{alg:quasigreedypayment} with the input $\Semi = \langle 3, 7, 10 \rangle$ and $k=28$}
    \label{tab:execution-algorithm-subtract} 
    \end{table}
\end{example}

\begin{definition}
    \label{def:quasigreedycost} 
    For a given numerical semigroup $\Semi = \langle S \rangle$ with generators $S=\{ s_1, s_2, \ldots, s_t \}$, where $1 \leq s_2 < \ldots < s_t$, and a given $k \in \langle S \rangle, k > 0$, the \emph{generalized greedy representation} of $k$ with respect to $S$, henceforth denoted $\mbox{\QGR}_S(k)$, is the representation vector $\textbf{a} = (a_1, a_2, \ldots, a_t)$ produced by Algorithm \ref{alg:quasigreedypayment}, and $\displaystyle \mbox{\QGC}_S(k) = \sum_{i=1}^{t} a_i$. 
\end{definition} 

All representable numbers $k \in \langle S \rangle, \ k > 0$, have a generalized greedy representation. Moreover, if $k$ is greedy-representable, then the generalized greedy representation is just the greedy representation. Thus, for the sake of brevity, from now on we will use the term \lq greedy\rq \ to mean \lq generalized greedy\rq. 

Notice also that if $s_1 = 1$ then $\langle S \rangle = \N_0$, and any $k \in \N$ is greedy-representable, hence the  greedy representation of $k$ can be obtained with the aid of Algorithm \ref{alg:quasigreedypayment}. Nonetheless, in this paper we will mostly be concerned with the case $s_1 > 1$.  

The drawback of Algorithm \ref{alg:quasigreedypayment} is perhaps its efficiency: it uses too many steps for obtaining the greedy representation of a number $k$. We can speed up the process by using division instead of subtraction, as shown in Algorithm \ref{alg:quasigreedypayment2}. \footnote{Someone may argue that division is more costly than subtraction, but we will defer these considerations for later.}
\\\\ 
	\begin{algorithm}[H] 
		\SetKwInOut{Input}{Input}
		\SetKwInOut{Output}{Output}
		\SetKw{DownTo}{downto}
		\vspace{.2cm}
		\Input{The set of denominations $S = \{ s_1, s_2, \ldots, s_t \}$, with $1 \leq s_1 < s_2 < \ldots < s_t$, $\gcd(s_1, s_2, \ldots, s_t) = 1$, and an element $k \in \langle S \rangle$, $k > 0$.} 
		\Output{Greedy representation vector $\textbf{a} = (a_1, a_2, \ldots, a_t)$.} 
		\vspace{.2cm}
		
		\For{ $i$:= $t$ \DownTo $1$} 
		{ 
            Let $q$ be the largest integer such that $k = q s_i + r$ and $r \in \langle S \rangle$ \; 
			 $a_i$ := $q$\;
			 $k$ :=  $r$\;  
              \If{ $k=0$ }{ \Return{ $\textbf{a}$ }\; } 
		} 
		\caption{GENERALIZED GREEDY REPRESENTATION METHOD (DIVISION VERSION)}
		\label{alg:quasigreedypayment2}  
	\end{algorithm}

\vspace{3mm}
\begin{example}
    \label{ex:algorithmic0}
    Let $S = \{ 3, 7, 10 \}$ be a set of denominations and $\Semi = \langle S \rangle$. If we apply Algorithm \ref{alg:quasigreedypayment2} to $k=28$ we observe the following steps: 
    \begin{table}[h]
    \centering
        \begin{tabular}{|c|c|c|c|c|}
            \hline
            Index $i$ & Target $k$ & Generator $s_i$ & Quotient $q$ & Remainder $r$ \\ \hline
            3  & 28  & 10 & 1  & 18   \\ 
            2  & 18  & 7  & 0  & 18  \\ 
            1  & 18  & 3  & 6  & 0  \\ \hline
        \end{tabular} 
    \caption{Execution of Algorithm \ref{alg:quasigreedypayment2} with the input $\Semi = \langle 3, 7, 10 \rangle$ and $k=28$}
    \label{tab:execution-algorithm3}
    \end{table}
    Notice that on the first line of the table (i.e. the one corresponding to $i=3$) the quotient is $q=1$, and not $q=2$, as we would expect in ordinary division, since if we take $q=2$ we would get a remainder $r=8$, which does not belong to $\langle S \rangle$. Similarly, on the second line we have $q=0$ because $q=1$ or $q=2$ would result in a remainder $r=11$ or $r=4$, respectively, none of them belonging to $\langle S \rangle$. Thus, the output of the algorithm is the representation vector $(6,0,1)$.     
\end{example}

Regardless of the algorithm we use for obtaining the greedy representation of a number, there is a straightforward observation about the greedy representation, which we will state in the form of a lemma, since it will become useful later. 

\begin{lemma}
    \label{lemma:singlerep}
    Let $\Semi = \langle S \rangle$ be a numerical semigroup with generating set $S=\{ s_1, s_2, \ldots, s_t \}$, such that $1 \leq s_1 < s_2 < \ldots < s_t$, and suppose also that $k \in \langle S \rangle, k > 0$. If $\mbox{\QGR}_S(k)$ involves only one $s_i$, then 
    \begin{displaymath}
        \mbox{\QGR}_S(k) = \mbox{\MREP}_S(k). 
    \end{displaymath}
\end{lemma}

\EndProof 

\vspace{2mm}
Let us now investigate the greedy representation a bit further. Given a representation vector $\textbf{a} = (a_1, a_2, \ldots, a_t)$, let us denote by $\hat{\textbf{a}}$ the \emph{reverse} of $\textbf{a}$, i.e. $\hat{\textbf{a}} = (a_t, a_{t-1}, \ldots, a_1)$. If $\textbf{a} = (a_1, a_2, \ldots, a_t)$ is the greedy representation of $k$, then $\hat{\textbf{a}}$ is the largest vector in lexicographic order among all vectors $\hat{\textbf{b}} = (b_t, b_{t-1}, \ldots, b_1)$ that represent $k$. The reverse representation vector $\hat{\textbf{a}}$ is just another way of looking at the same data. Sometimes it is more convenient to work with $\hat{\textbf{a}}$, rather than working with $\textbf{a}$, as in \cite{Pea05}, and in Proposition \ref{prop:cardinalitytwo} below.  

Obviously, the greedy representation is not necessarily the best or the most efficient representation of $k$, as we saw in Example \ref{ex:algorithmic-subtract}. That is, given the greedy representation vector $\textbf{a} = (a_1, a_2, \ldots, a_t)$, the sum $\displaystyle \sum_{i=1}^{t} a_i$ need not be minimal. However, sometimes it is, which leads us to the following concept:   

\begin{definition}
    \label{def:quasigreedyset} 
    Let $\Semi = \langle S \rangle$ be a numerical semigroup with generating set $S = \{ s_1, s_2, \ldots, s_t \}$, where $1 \leq s_1 < s_2 < \ldots < s_t$, such that Algorithm \ref{alg:quasigreedypayment} \emph{always} produces an optimal representation for \emph{any} given $k \in \Semi$. Then $S$ will be  called \emph{ greedy (in the general sense)}, and the semigroup $\Semi = \langle S \rangle$ will also be called \emph{ greedy (in the general sense)}. 
\end{definition}
Sets of cardinality two provide the first examples of greedy sets: 
\begin{proposition}
    \label{prop:cardinalitytwo}
    Let $S = \{ s_1, s_2 \}$, with $1 < s_1 < s_2$ and $\gcd(s_1, s_2) = 1$. Then $S$ is greedy. 
\end{proposition} 

\noindent \textbf{Proof:} Let $k$ be an arbitrary element of $\langle S \rangle$, with $k > 0$, and let $\hat{\textbf{a}} = (a_2, a_1)$ be the reverse greedy representation of $k$. Moreover, let $\hat{\textbf{b}} = (b_2, b_1)$ be another representation of $k$, and suppose that $\hat{\textbf{b}}$ is more efficient than $\hat{\textbf{a}}$, i.e. $b_2 + b_1 < a_2 + a_1$. Since $\hat{\textbf{a}}$ is the largest representation in lexicographic order, then $a_2 > b_2$ which means that $b_1 > a_1$. Thus, from $b_2 + b_1 < a_2 + a_1$ we get that $0 < b_1 - a_1 < a_2 - b_2$. 

On the other hand, since both $\hat{\textbf{a}}$ and $\hat{\textbf{b}}$ are representations of $k$ we have 
\begin{align*}
    k &= a_2 s_2 + a_1 s_1 \\
    k &= b_2 s_2 + b_1 s_1. 
\end{align*}
Subtracting the second equation from the first we get 
\begin{displaymath}
    (a_2 - b_2) s_2 = (b_1 - a_1) s_1. 
\end{displaymath}
Since $s_2 > s_1$, then $a_2 - b_2 < b_1 - a_1$, which is a contradiction. 

\EndProof 

Proposition \ref{prop:cardinalitytwo} answers the question of existence, and it raises naturally further questions. To begin with, we want to know what happens in the case of three (or more) generators. In this case we can no longer expect the system of generators to be automatically greedy, as the next example shows: 
\begin{example}[Non-greedy sets of cardinality 3]
\label{example:non-quasi-greedy1} 
Take $S = \{ 5, 8, 9 \}$ and $k=24$. It can be verified that $\mbox{\MREP}_S(24) = (0,3,0)$, whereas on the other hand $\mbox{\QGR}_S(24) = (3,0,1)$. So, $\mbox{\MC}_S(24) = 3 < \mbox{\QGC}_S(24) = 4$. Another non-greedy set of cardinality three is the set $S = \{ 3, 7, 10 \}$ of Example \ref{ex:algorithmic-subtract}.  
\end{example}

\section{Counterexamples, witnesses and algorithmic identification of greedy sets} 
\label{sec:counterwitnesses} 

If a set $S$ is not greedy (in the general sense), then there must exist some $k$ such that $\mbox{\MC}_S(k) < \mbox{\QGC}_S(k)$. We call any such $k$ a \emph{counterexample}. Note that, by Lemma \ref{lemma:singlerep}, if $k$ is a counterexample then $\mbox{\QGR}_S(k)$ must involve at least two generators. 

The second observation, which we will now investigate in more detail, is that the smallest counterexample must lie in some finite interval, which we call the \emph{critical range}. This will be the basis for the algorithmic identification of greedy sets: If we cannot find any counterexample in the critical  range then we can conclude that $S$ is greedy. 

The following result is a straightforward generalization of Theorem 2.2 of \cite{KoZaks94} to the case of numerical semigroups: 
\begin{theorem}
    \label{theo:range}
    Let $\Semi = \langle S \rangle$ be a numerical semigroup generated by $S = \{ s_1, s_2, \ldots, s_t \}$, with $1 < s_1 < s_2 < \cdots < s_t$. If there exists a counterexample $k \in \Semi$ such that $\mbox{\MC}_S(k) < \mbox{\QGC}_S(k)$, then the smallest such $k$ lies in the range 
    \begin{equation}
        \label{eq:range}
        s_3 + s_1 + 2 \leq k \leq F(\Semi) + s_t + s_{t-1}. 
    \end{equation}
    Moreover, these bounds are tight. 
\end{theorem} 
\vspace{4mm} 
\noindent \textbf{Proof:} On one hand, according to Proposition \ref{prop:cardinalitytwo}, if $k < s_3$ is representable, then $\mbox{\MC}_S(k) = \mbox{\QGC}_S(k)$. Similarly, $\mbox{\MC}_S(s_3) = \mbox{\QGC}_S(s_3) = 1$. On the other hand, any representable $k \in \{ s_3+1, \ldots, s_3+s_1-1 \}$ cannot be a counterexample, since $k-s_3$ is not representable, and therefore any representation of $k$ would only use the generators $\{ s_1, s_2 \}$, which is a greedy set. Similarly, $s_3 + s_1$ is not a counterexample because $\mbox{\MC}_S(s_3 + s_1) = \mbox{\QGC}_S(s_3 + s_1) = 2$. 

Finally, let us verify that $k = s_3 + s_1 + 1$ cannot be a counterexample either. Let us assume that $S = \{ s_1, s_2, s_3 \}$. The only way that the generalized greedy representation of $k$ involves $s_3$ is by letting $s_2 = s_1 + 1$, but in that case $\mbox{\MC}_S(s_3 + s_1 + 1) = \mbox{\QGC}_S(s_3 + s_1 + 1) = 2$. If $s_2 > s_1 + 1$, then $s_3$ will not appear in the greedy representation of $k$, and $k$ cannot be a counterexample because $\{ s_1, s_2 \}$ is a greedy set. 

Now suppose that $S = \{ s_1, s_2, s_3, s_4 \}$, with $s_3 + 1 \leq s_4 \leq k$. If $s_4 = s_3 + 1$ we have $\mbox{\MC}_S(s_3 + s_1 + 1) = \mbox{\QGC}_S(s_3 + s_1 + 1) = 2$. In the other extreme, if $s_4 = k = s_3 + s_1 + 1$ we have $\mbox{\MC}_S(s_3 + s_1 + 1) = \mbox{\QGC}_S(s_3 + s_1 + 1) = 1$. So, let $s_3 + 1 <  s_4 < s_3 + s_1 + 1$. Then we have $0 <  k-s_4 < s_1$, which means that $k-s_4 \notin \Semi$, and thus $s_4$ cannot appear in the greedy representation of $s_3 + s_1 + 1$, which takes us  back to the case $S = \{ s_1, s_2, s_3 \}$. 

The same reasoning applies if we have more than $4$ generators. This establishes the lower bound. 

As for the upper bound, let $k \in \Semi, \ k > F(\Semi) + s_t + s_{t-1}$ and assume inductively that $\mbox{\MC}_S(j) = \mbox{\QGC}_S(j)$ for all representable $j < k$. Now let $s_i$ be \emph{any} generator used in some \emph{optimal} representation of $k$. If $i=t$ then
\begin{align*}
    \mbox{\QGC}_S(k) &= \mbox{\QGC}_S(k-s_t) + 1 \\ 
                    & \; \mbox{ (by definition of \GRC) } \\ 
                    &= \mbox{\MC}_S(k-s_t) + 1 \\
                    & \; \mbox{ (by the induction hypothesis) } \\
                    &= \mbox{\MC}_S(k) \; \mbox{ (by Lemma \ref{lemma:kozaks}) } 
\end{align*}


On the other hand, if $i<t$ then 

\begin{align*}
    \mbox{\QGC}_S(k) &= \mbox{\QGC}_S(k-s_t) + 1 \\ 
                    & \quad \mbox{ (by definition of \QGC) } \\ 
                    &= \mbox{\MC}_S(k-s_t) + 1 \\
                    & \quad \mbox{ (by the induction hypothesis) } \\
                    &\leq \mbox{\MC}_S(k-s_t-s_i) + 2 \\
                    & \quad \mbox{ (by Lemma \ref{lemma:kozaks}) } \\
                    &\leq \mbox{\QGC}_S(k-s_t-s_i) + 2 \\
                    & \quad \mbox{ (by definition of \MC) } \\ 
                    &= \mbox{\QGC}_S(k-s_i) + 1 \\
                    & \quad \mbox{ (by definition of \QGC) } \\ 
                    &= \mbox{\MC}_S(k-s_i) + 1 \\
                    & \quad \mbox{ (by the induction hypothesis) } \\
                    &= \mbox{\MC}_S(k) \\
                    & \quad \mbox{ (by Lemma \ref{lemma:kozaks}) } \\ 
                    &\leq \mbox{\QGC}_S(k) \\
                    & \quad \mbox{ (by definition of \MC). } 
\end{align*} 
Thus, in either case $\mbox{\MC}_S(k) = \mbox{\QGC}_S(k)$. 
\\\\ 
Now we have to show that these bounds are tight. Let us start with the lower bound, which is easier. Consider $S = \{ 2, j, 2j-4 \}$, with $j \geq 7$ odd. The lower bound for this system is $s_3 + s_1 + 2 = 2j$. Now check that $\mbox{\QGC}_S(2j) = 3$, whereas $\mbox{\MC}_S(2j) = 2$. Therefore, $2j$ is a counterexample. 

In contrast, proving that the upper bound is tight is much more involved. Consider $S = \{ 2, j, j+1 \}$, with $j \geq 6$ even. The Frobenius number of $\langle S \rangle$ is $j-1$, and the critical range is $[j+5; 3j]$. It is easy to see that $\mbox{\MC}_S(3j) = 3$. On the other hand, $\mbox{\QGC}_S(3j) = \frac{j}{2} + 1 \geq 4$. Therefore, the upper bound $3j$ is a counterexample. Now we have to show that it is the \emph{smallest} counterexample. 

Let's break the interval $[j+5; 3j-1]$ into two sub-intervals or cases:
\begin{enumerate}[(A)]   
    \item \label{item:first} $[j+5; 2j]$, and
    \item \label{item:second} $[2j+1; 3j-1]$. 
\end{enumerate} 
The numbers of case \ref{item:first} do not have a representation involving the generator $j+1$, therefore, they belong to $\langle 2, j \rangle$, which is greedy (see Proposition \ref{prop:cardinalitytwo}). 

In case \ref{item:second} we can again distinguish several subcases: 
\begin{enumerate}[(B1)]
    \item \label{item:subfirst} $k=2j+1$:  $\mbox{\MC}_S(k) =  \mbox{\QGC}_S(k) = 2$ 
    \item \label{item:subsecond} $k=2j+2$: $\mbox{\MC}_S(k) =  \mbox{\QGC}_S(k) = 2$   
    \item \label{item:subthird} $2j+3 \leq k \leq 3j-1$, with $k$ odd: In this case $\displaystyle \mbox{\QGR}_S(k) = 1 \cdot (j+1) + 1 \cdot j + \left( \frac{k-2j-1}{2} \right) \cdot 2$, hence $\displaystyle \mbox{\QGC}_S(k) = \frac{k-2j-1}{2} + 2$. We can show that this is also $\mbox{\MC}_S(k)$. Indeed, since $k$ is odd, \emph{any} representation of $k$ must use the generator $j+1$ an odd number of times. However, $j+1$ cannot appear three  times, since $k < 3j$, so it can only appear once. Now, $k-(j+1) \in \langle 2, j \rangle$, which is greedy. 
    \item \label{item:subfourth} $2j+4 \leq k \leq 3j-2$, with $k$ even: $\displaystyle \mbox{\QGR}_S(k) = 2 \cdot (j+1) + \left( \frac{k-2j-2}{2} \right) \cdot 2$, hence $\displaystyle \mbox{\QGC}_S(k) = \frac{k-2j-2}{2} + 2$. We can also show that this is $\mbox{\MC}_S(k)$. Since $k$ is even, the optimal representation of $k$ must involve two occurrences of the generator $j$, or two occurrences of the generator $j+1$ (it cannot involve a combination of $j$ and $j+1$). Clearly, the representation that uses $j+1$ will be more efficient than the former. The rest is obvious. 
\end{enumerate} 
\EndProof 



\vspace{2mm} 
Theorem \ref{theo:range} is the starting point for the algorithmic identification of greedy sets. The prospective algorithm must look for a counterexample in the range provided by the inequalities (\ref{eq:range}), and if we cannot find one, then we can conclude that the given set $S$ is greedy. The next example illustrates the procedure. 

\begin{example}
    \label{ex:algorithmic1}
    Let $S_1 = \{ 3, 7, 10 \}$ and $S_2 = \{ 3, 7, 11 \}$ be two sets of generators of cardinality three, and $\Semi_1 = \langle S_1 \rangle$ and $\Semi_2 = \langle S_2 \rangle$ the numerical semigroups generated by $S_1$ and $S_2$ respectively. Their Frobenius numbers are $F(\Semi_1) = 11$ and $F(\Semi_2) = 8$, and their respective critical ranges are $R_1 = [15; 28]$ and $R_2 = [16; 26]$. 
    
    In the case of $\Semi_1$ we can verify that $\mbox{\MC}_{S_1}(k) = \mbox{\QGC}_{S_1}(k)$ for all $k \in [15; 27]$. For $k=28$ we saw in Example \ref{ex:algorithmic0} that the greedy representation vector is $(6,0,1)$, hence $\mbox{\QGC}_{S_1}(28) = 7$. However, $28 = 4 \cdot 7$, hence $\mbox{\MC}_{S_1}(28) = 4$. Thus, $28$ is a counterexample, and $\Semi_1$ is not greedy. 

    On the other hand, in the case of $\Semi_2$ we can verify that $\mbox{\MC}_{S_1}(k) = \mbox{\QGC}_{S_2}(k)$ for all $k$ in the critical range. Hence, $\Semi_2$ is greedy. 
    
\end{example}

There is a catch to Theorem \ref{theo:range}, though: Looking for a counterexample involves calculating the minimal representation of all $k$ in the range (\ref{eq:range}), which may be a costly process, as shown in \cite{Lue75}. Kozen and Zaks \cite{KoZaks94} circumvent the calculation of $\mbox{\MC}_S(k)$ by introducing the concept of a \emph{witness}, which we can also adapt to our setting with small changes.   
\begin{definition}
    \label{def:witness} 
    Given a set of generators $S=\{ s_1, s_2, \ldots, s_t \}$, with $1 < s_1 < s_2 < \ldots < s_t$ and $\gcd(s_1, s_2, \ldots, s_t) = 1$, a \emph{witness} for $S$ is any representable integer $k > 0$, such that
    \begin{displaymath}
        \mbox{\QGC}_S(k) > \mbox{\QGC}_S(k-s_i)+1 
    \end{displaymath}
     for some generator $s_i < k$.  
\end{definition}
\begin{lemma}  
    \label{lemma:witness}
    As in Definition \ref{def:witness}, let $S=\{ s_1, \ldots, s_t \}$ be a set of nonnegative integers, such that $1 < s_2 < \ldots < s_t$ and $\gcd(s_1, \ldots, s_t) = 1$. Then 
    \begin{itemize} 
        \item Every witness for $S$ is a counterexample. 
        \item The smallest counterexample (if it exists) is also a witness. 
    \end{itemize}
\end{lemma}
\noindent \textbf{Proof:} 
\begin{itemize}
    \item In order to prove the first claim, suppose that $k$ is a witness for $S$. Then $\mbox{\QGC}_S(k) > \mbox{\QGC}_S(k-s_i)+1$ for some generator $s_i$. Note that this implies that $k-s_i$ is representable.  Therefore, 
    \begin{align*}
        \mbox{\MC}_S(k) &\leq \mbox{\MC}_S(k-s_i)+1 \quad \mbox{ (by Lemma \ref{lemma:kozaks}) } \\
                        &\leq \mbox{\QGC}_S(k-s_i)+1 \; \mbox{ (by def. \MC) } \\
                        &< \mbox{\QGC}_S(k) \quad \mbox{ (by definition of witness). } 
    \end{align*} 
    \item For the second claim, let $k$ be a counterexample that is not a witness, and let $s_i$ be any generator used in some optimal representation of $k$. We will now show that $k-s_i$ is also a counterexample: 
    \begin{align*}
        \mbox{\MC}_S(k-s_i) &= \mbox{\MC}_S(k)-1 \quad \quad \quad \mbox{ (by Lemma \ref{lemma:kozaks}) } \\
                        &< \mbox{\QGC}_S(k)-1  \; \mbox{ (by def. \MC } \\
                        &  \quad \quad \quad \quad \quad \quad \mbox{ and the fact that $k$ is counterexample) } \\ 
                        &\leq \mbox{\QGC}_S(k-s_i),       
    \end{align*}  
    i.e. $\mbox{\MC}_S(k-s_i) < \mbox{\QGC}_S(k-s_i)$. Therefore, the smallest counterexample must also be a witness. 
\end{itemize}

\EndProof 
\\
An immediate consequence of Theorem \ref{theo:range} and Lemma \ref{lemma:witness} is the following
\begin{theorem}
    \label{theo:witness}
    As in Definition \ref{def:witness}, let $S=\{ s_1, \ldots, s_t \}$, with $1 < s_1 < s_2 < \ldots < s_t$ and $\gcd(s_1, \ldots, s_t) = 1$, and let $\Semi = \langle S \rangle$. Then, $S$ is greedy if, and only if, $S$ does not have any witness $k$ in the interval 
    \begin{displaymath}
        s_3 + s_1 + 2 \leq k \leq F(\Semi) + s_t + s_{t-1}.  
    \end{displaymath}    
\end{theorem} 
\EndProof 

\section{Complexity issues}
\label{sec:complexity}

Theorem \ref{theo:witness} finally provides a practical algorithm for checking the greediness of some set of generators $S$. Now there is no need to compute the optimal representation for any $k$. The algorithm just verifies whether the inequality 
\begin{displaymath}
    \mbox{\QGC}_S(k) \leq \mbox{\QGC}_S(k-s_i)+1
\end{displaymath} 
holds for all $k$ in the range $s_3 + s_1 + 1 < k \leq F(\Semi) + s_t + s_{t-1}$ and all generators $s_i < k$. These ideas are formalized in Algorithm \ref{alg:search-witness}. 

	\begin{algorithm}[ht] 
		\SetKwInOut{Input}{Input}
		\SetKwInOut{Output}{Output}
		\SetKw{DownTo}{downto}
        \SetKw{To}{to} 
		\vspace{.2cm}
		\Input{The set of denominations $S = \{ s_1, s_2, \ldots, s_t \}$, with $1 < s_1 < s_2 < \ldots < s_t$, $\gcd(s_1, s_2, \ldots, s_t) = 1$.}  
		\Output{TRUE if $\langle S \rangle$ is greedy, and  FALSE otherwise.} 
		\vspace{.2cm}

        $\Semi$ := $\langle S \rangle$\;  
		\For{ $k$:= $s_3 + s_1 + 2$ \To $F(\Semi) + s_t + s_{t-1}$} 
		{ 
			 $t'$ :=  Largest index $j$ such that $s_{j} \leq k$\;  
                \For{ $i$:= $1$ \To $t'-1$}  
                { 
                    \If{ $\mbox{\QGC}_S(k) > \mbox{\QGC}_S(k-s_i)+1$ }{ \Return{ FALSE }\; } 
                }
                
		} 
        \Return{ TRUE }\;
		\caption{DETERMINE WHETHER A SEMIGROUP DEFINED BY SET OF GENERATORS IS GREEDY}
		\label{alg:search-witness}  
	\end{algorithm}

\begin{observation}
    Notice that the inner loop of Algorithm \ref{alg:search-witness} runs through the indices $1, \ldots, t'-1$, instead of $1, \ldots, t'$. Although the definition of witness requires that 
    \begin{displaymath}
        \mbox{\QGC}_S(k) > \mbox{\QGC}_S(k-s_i)+1 
    \end{displaymath} 
    for \emph{some} generator $s_i$, in practice it does not make sense to test the largest possible generator $s_i < k$ (or $s_i = k$, for that matter), since in that case 
    \begin{displaymath}
        \mbox{\QGC}_S(k) = \mbox{\QGC}_S(k-s_i)+1.  
    \end{displaymath}
    Subtracting the largest possible generator is precisely what the greedy algorithm does. 
\end{observation}

\begin{example}
    \label{ex:algorithmic2}
    In Example \ref{ex:algorithmic1} we investigated the numerical semigroups $\Semi_1 = \langle 3, 7, 10 \rangle$ and $\Semi_2 = \langle 3, 7, 11 \rangle$, with respective Frobenius numbers $F(\Semi_1) = 11$ and $F(\Semi_2) = 8$, and respective critical ranges $R_1 = [15; 28]$ and $R_2 = [16; 26]$. As we saw, $\mbox{\QGC}_{S_1}(28) = 7$ and $\mbox{\QGR}_{S_1}(28) = (6, 0, 1)$. However, if we now apply Algorithm \ref{alg:search-witness} on $\Semi_1$ we can see that $\mbox{\QGR}_{S_1}(28-7) = \mbox{\QGR}_{S_1}(21) = (0, 3, 0)$. Therefore, $\mbox{\QGC}_{S_1}(21) = 3$ and 
     \begin{displaymath}
        \mbox{\QGC}_{S_1}(28) > \mbox{\QGC}_{S_1}(28-7)+1 = 4.  
    \end{displaymath}
    Hence $k=28$ is a witness, and we can conclude that $\Semi_1$ is not greedy. 

    \noindent On the other hand, if we apply Algorithm \ref{alg:search-witness} on $\Semi_2$, we can verify that there are no witnesses in the critical range $R_2$, hence $\Semi_2$ is greedy. 
\end{example}

In the worst case, the outer loop of Algorithm \ref{alg:search-witness} runs through the whole critical range, whose number of elements, $N = F(\Semi) + s_t + s_{t-1} - s_1 - s_3 - 1$, is polynomial in $s_t$, the largest generator. On the other hand, the inner loop runs through all the $t$ generators in the worst case. Finally, the calculation of $\mbox{\QGC}_S(k)$ also runs in polynomial time. Therefore, the complexity of Algorithm \ref{alg:search-witness} is polynomial in $t s_t$. However, this running time is potentially exponential in the number of digits of the input, if the input is represented in binary or decimal. 

In the original change-making problem, i.e. when $s_1 = 1$, there does exist an algorithm with a running time that is polynomial in the number of digits \cite{Pea05}, but we still don't know if that algorithm (or a modified version) is applicable in our case.  

Regarding the complexity of other related questions, Lueker proved that finding an optimal representation of $k$ is NP-hard in the original change-making problem, i.e. when $s_1 = 1$ (its decision version is NP-complete) \cite{Lue75}. Therefore, it must also be NP-hard in our more general setting, when $s_1 \geq 1$. 

On the other hand, Kozen and Zaks proved that determining whether the greedy representation of $k$ is optimal is co-NP-complete in the original setting \cite{KoZaks94}. Note that the greedy representation of $k$ is optimal if there does not exist any representation of $k$ that is more efficient than the greedy one. The complementary problem is to determine if the greedy representation of $k$ is \emph{not} optimal, which means that there exists a better representation of $k$ than the greedy one. That amounts to finding the optimal representation of $k$, which is NP-hard, as we have just said above. 

The same reasoning applies to our generalized setting, i.e. $s_1 \geq 1$. That is,  the question of deciding whether the greedy representation of $k$ is optimal is co-NP-complete. 

\section{Semigroups of embedding dimension three}
\label{sec:threegens} 

We have already seen that semigroups generated by two generators are greedy. For three generators this is not always the case. We will now investigate the structure of potential counterexamples in the case of three generators, and thus find specific criteria to decide greediness in that particular case. 

\begin{lemma}
    \label{lemma:threegens} 
    Let $S = \{ s_1, s_2, s_3 \}$, with $1 < s_1 < s_2 < s_3$, and $\gcd(s_1, s_2, s_3) = 1$, so that $\Semi = \langle S \rangle$ is not greedy. Then, the smallest counterexample $k$ (and hence the smallest witness) has the form $k = s_2 y$, and it is a  solution of the Diophantine equation 
    \begin{equation}
        s_1 x + s_3 z = s_2 y, 
    \end{equation}
    where $x, y, z$ are positive integers, such that $y < x+z$. 
\end{lemma}  

\noindent \textbf{Proof:} \\ 
Let $k$ be the smallest counterexample of $\Semi = \langle S \rangle$ (and hence also the smallest witness). Suppose that $\mbox{\QGR}_S(k) = \alpha s_1 + \beta s_2 + \gamma s_3$, where $\gamma > 0$, for otherwise $\mbox{\QGR}_S(k) = \alpha s_1 + \beta s_2$ would be the optimal representation of $k$, and hence $k$ would not be a counterexample. 

Now let $\mbox{\MC}_S(k) = \alpha' s_1 + \beta' s_2 + \gamma' s_3$, so that $\alpha'+\beta'+\gamma' < \alpha+\beta+\gamma$. We must have $\delta' < \delta$ for otherwise $\alpha' s_1 + \beta' s_2 + \gamma' s_3$ would be the $\mbox{\QGR}_S(k)$. So let us assume that $0 < \delta' < \delta$. We can now subtract $\delta' s_3$ from both representations to obtain
$k-\delta' s_3 = \alpha s_1 + \beta s_2 + (\gamma-\gamma') s_3$ on one side and $k-\delta' s_3 = \alpha' s_1 + \beta' s_2$ on the other side, where $\alpha' + \beta' < \alpha + \beta + (\gamma-\gamma')$ and $\gamma-\gamma' > 0$. Consequently, $k-\delta' s_3$ would also be a counterexample, which contradicts the assumption that $k$ was the smallest counterexample. Therefore we must conclude that $\delta' = 0$. 

So now we have two representations of $k$: $k = \alpha s_1 + \beta s_2 + \gamma s_3$ and $k = \alpha' s_1 + \beta' s_2$, where $\alpha' + \beta' < \alpha + \beta + \gamma$. For that to happen we must have $\beta < \beta'$. Let's suppose that $0 < \beta < \beta'$. In that case we may subtract $\beta s_2$ from both representations to obtain $k-\beta s_2 = \alpha s_1 + \gamma s_3$ on one side, and $k-\beta s_2 = \alpha' s_1 + (\beta' - \beta) s_2$ on the other side, where $\alpha' + \beta' - \beta < \alpha + \gamma$. That means $k-\beta s_2$ would also be a counterexample, which contradicts the assumption that $k$ is the smallest counterexample. That leads us to conclude that $\beta = 0$. 

So, both representations of $k$ now are: $k = \alpha s_1 +  \gamma s_3$ and $k = \alpha' s_1 + \beta' s_2$, where $\alpha' + \beta' < \alpha + \gamma$. For this to happen we must have $\alpha' < \alpha$ and $\beta' > \gamma$. Let's now consider the number $k-\alpha' s_1$. On one side $k-\alpha' s_1 = \beta' s_2$, and on the other side $k-\alpha' s_1 = (\alpha-\alpha') s_1 + \gamma s_3$, where $\beta' < \alpha-\alpha'+\delta$. That means $k-\alpha' s_1$ would also be a counterexample, which again contradicts the fact that $k$ is the smallest counterexample. Hence $\alpha' = 0$. 

In the end we have $k = \alpha s_1 +  \gamma s_3$ and $k = \beta' s_2$, where $\beta' < \alpha + \gamma$ (since $k$ is a counterexample). Equate both expressions of $k$, replacing $\alpha$ by $x$, $\beta'$ by $y$ and $\gamma$ by $z$, and we get our Diophantine equation. 

\EndProof

This lemma suggests a simpler algorithm for deciding greediness in the case of three generators. We just have to look for all the multiples $j s_2$ of $s_2$ that lie in the critical range, and check whether the greedy cost of $j s_2$ is \emph{worse} than $j$ (i.e. $>j$). Then, this $k = j s_2$ would be our counterexample. If we cannot find such a $k$, then the semigroup is greedy. 

We can now apply the above results to particular families of semigroups with three generators. For instance, 

\begin{theorem}
    \label{theo:threeconsecutivegens} 
    Let  $\Semi = \langle n, n+1, n+2 \rangle$. Then, $\Semi$ is greedy. 
\end{theorem}

\noindent \textbf{Proof:} \\  
According to Lemma \ref{lemma:threegens}, any potential counterexample $k$ must have the form $y(n+1)$, and it must satisfy 
\begin{displaymath}
    x n + z (n+2) = (x+z)n + 2z = y(n+1) = yn + y.
\end{displaymath} 
We can easily find a subset of the solutions by decomposing the equation into two, namely $y = x+z$ and $y=2z$. The first equation is already telling us that $k$ is not a counterexample. Additionally, this implies that $y$ must be even, and $\displaystyle x=z=\frac{y}{2}$. Let us now check for other potential solutions: 
\begin{align*}
    y(n+1) &= x n + z (n+2) = n(x+z) + 2z \\
            &= n(x+z) + (x+z) + 2z - (x+z) \\
            &= (n+1)(x+z) + z-x 
\end{align*}
This implies that $z-x$ must be a multiple of $n+1$, i.e. $z-x = m(n+1)$, where $m \in \Z$ (not necessarily positive). Notice that if we take $m=0$ we will just obtain the first subset of solutions above, which as we saw, do not yield any counterexamples. 

Let us now consider two cases: positive $m$ and negative $m$. We begin with positive $m$, which implies that $z>x$. More precisely,
\begin{align*}
    z &= x + m(n+1) \\
    y &= 2x + m(n+2). 
\end{align*} 
Recall that any potential counterexamples $k=y(n+1)$ must lie in the critical range. The Frobenius number for $\Semi = \langle n, n+1, n+2 \rangle$ is $\displaystyle \Bigl \lfloor \frac{n}{2} \Bigr \rfloor n - 1$ (see \cite{Kan97}, for instance), hence the critical range for $\Semi = \langle n, n+1, n+2 \rangle$ is 
\begin{displaymath}
    \biggl[ 2n+4; \Bigl \lfloor \frac{n}{2} \Bigr \rfloor n + 2n + 2  \biggr] 
\end{displaymath}
So, in particular 
\begin{align*}
    y(n+1) \, \leq \, \Bigl \lfloor \frac{n}{2} \Bigr \rfloor n + 2n + 2 \, &\leq \, \frac{1}{2} (n+2)^2 \\ 
    (n+1)\left( 2x + m(n+2) \right) \, &\leq \, \frac{1}{2} (n+2)^2, 
\end{align*}
which leads to 
\begin{align*}
    x \, &\leq \, \frac{(n+2)^2-2m(n+1)(n+2)}{4(n+1)} \\
    x \, &\leq \, \frac{-2mn^2-6mn-4m+n^2+4n+4}{4(n+1)} 
\end{align*}
It is not difficult to verify that the right hand side of this inequality takes value $1-m$ at $n=0$ and it is decreasing as $n$ increases, which means that the variable $x$ must take negative values if $y(n+1)$ were to lie inside the critical range. Hence, these solutions do not lead to any counterexamples. 

On the other hand, if $m$ is negative we get 
\begin{align*}
    x &= z + m'(n+1) \\
    y &= 2x + m'n, 
\end{align*}
where $m'=-m$ is positive, which leads to 
\begin{align*}
    y(n+1) \, &\leq \, \frac{1}{2} (n+2)^2 \\ 
    (n+1)\left( 2z + m'n \right) \, &\leq \, \frac{1}{2} (n+2)^2,  
\end{align*}     
whence 
\begin{align*}    
    z \, &\leq \, \frac{(n+2)^2-2m'n(n+1)}{4(n+1)} \\
    z \, &\leq \, \frac{-2m'n^2-2m'n+n^2+4n+4}{4(n+1)}.
\end{align*}
The right-hand side of the inequality takes value $1$ at $n=0$ and it is decreasing as $n$ increases. Hence, in order that $y(n+1)$ lies inside the critical range, the variable $z$ must take negative values (or $0$) for $n>0$. In other words, these solutions do not lead to any counterexamples either. 

\EndProof 



\section{Some open problems}
\label{sec:open}

In this section we collect several interesting open questions that have arisen throughout our research. 

Starting with the algorithmic viewpoint, we would like to determine whether a semigroup is greedy with the highest possible efficiency. As mentioned in Section \ref{sec:complexity}, the running time of Algorithm \ref{alg:search-witness} may still be exponential in the number of digits if the numbers are represented in binary or decimal (or in any other base larger than $1$). For the original setting of the change-making problem there does exist an algorithm that runs in polynomial time \cite{Pea05}. The question here is whether that polynomial-time algorithm can be adapted to our setting without a significant loss in efficiency. 

Greedy semigroups seem to be relatively abundant. With the aid of Algorithm \ref{alg:search-witness} we have sampled $90$ semigroups generated by three small elements (between $2$ and $15$), and out of these, only $25$ of them were \emph{not} quasi-greedy. On the other hand, proving that a given infinite family of semigroups is greedy seems to be a much greater challenge. In Section \ref{sec:threegens} we have proved that numerical semigroups generated by three consecutive integers are greedy, but the problem is still open for other families.  

In fact, there are many other families of numerical semigroups whose Frobenius number is known, and hence their critical range, which is a starting point for investigating their greediness. Only in embedding dimension three, for instance, we have several families of semigroups whose Frobenius number is known, such as semigroups generated by Fibonacci or Lucas triples \cite{Fel09}, generated by three consecutive squares or cubes \cite{Le15}, by sexy prime triplets \cite{Hwa23}, semigroups where the three generators are pairwise relatively prime \cite{RoRo12}, and others \cite{Kan97}.  

Similarly, we want to know if the property of greediness is inherited under some semigroup operations, such as duplication \cite{Danna13}, shifting \cite{Cona18} and quotient \cite{Stra15, Bo24}. Finally, we want to establish the connection between  greediness and other semigroup invariants that are also related to the factorization of semigroup elements (e.g. the elasticity, Betti elements, tame degree and catenary degree \cite{Chap18, Chap21}).

\section*{Acknowledgements} 

The authors have been partially supported by Grant 2021 SGR 00115 from the Government of Catalonia, by Project ACITECH PID2021-124928NB-I00, funded by MCIN/AEI/ 10.13039/501100011033/FEDER, EU, and by Project HERMES, funded by the European Union NextGenerationEU/PRTR via INCIBE. 




\begin{thebibliography}{00}

\bibitem{Ada10} Adamaszek,~A. and M.~Adamaszek: \lq\lq Combinatorics of the change-making problem\rq\rq. {\it European Journal of Combinatorics} \textbf{31}, 47--63 (2010). 





\bibitem{Bo24} Bogart,~T., C.~O'Neill and K.~Woods: When is a numerical semigroup a quotient? {\it Bulletin of the Australian Mathematical Society} 109, 67--76 (2024). DOI: 10.1017/S0004972723000035. 

\bibitem{Cai09} Cai,~X.: \lq\lq Canonical Coin Systems for Change-Making Problems\rq\rq. {\it Procs. of the 9th IEEE Int. Conf. on Hybrid Intelligent Systems}, 499--504 (2009).


\bibitem{Chap18} Chapman,~S. and C.~O'Neill: Factoring in the Chicken McNugget Monoid. {\it Mathematics Magazine}, \textbf{91}(5) 323--336 (2018).

\bibitem{Chap21} Chapman,~S., P.~Garc\'{\i}a-S\'anchez and C.~O'Neill: Factoring in the Chicken McNugget Monoid. {\it Mathematics Magazine}, \textbf{91}(5) 323--336 (2018).


\bibitem{Cona18} Conaway,~R., F.~Gotti, J.~Horton, C.~O'Neill, R.~Pelayo, M.~Pracht and B.~Wissman: \lq\lq Minimal presentations of shifted numerical monoids\rq\rq. {\it International Journal of Algebra and Computation} \textbf{28}(1) (2018), 53--68. DOI: 10.1142/S0218196718500030. 

\bibitem{Cow08} Cowen,~L.J., R.~Cowen and A.~Steinberg: \lq\lq Totally Greedy Coin Sets and Greedy Obstructions\rq\rq. {\it The Electronic Journal of Combinatorics} \textbf{15} (2008), $\#$R90.








\bibitem{Danna13} D'Anna,~M. and F.~Strazzanti: The numerical duplication of a numerical semigroup. {\it Semigroup Forum} 87, 149--160 (2013). DOI: 10.1007/s00233-012-9451-x.










\bibitem{Fel09} Fel,~L.: Symmetric Numerical Semigroups Generated by Fibonacci and Lucas Triples. {\it Integers} 9, 107--116 (2009). DOI: 10.1515/ INTEG.2009.010














\bibitem{Hwa23} Hwang,~W.T. and K.~Song: The Frobenius problems for Sexy Prime Triplets. {\it International Journal of Mathematics and Computer Science} 18(2), 321--333 (2023). 


\bibitem{Kan97} Kan,~I.D., B.S.~Stechkin and I.V.~Sharkov: Frobenius Problem for Three Arguments. {\it Mathematical Notes} 62(4), 521--523 (1997). 

\bibitem{KoZaks94} Kozen,~D. and S.~Zaks: Optimal bounds for the change-making problem. {\it Theoretical Computer Science} 123, 377--388 (1994). 

\bibitem{Le15} Lepilov,~D., J.~O'Rourke and I.~Swanson: Frobenius numbers of numerical semigroups generated by three consecutive squares or cubes. {\it Semigroup Forum} 91, 238--259 (2015). DOI: 10.1007/s00233-014-9687-8. 




\bibitem{Lue75} Lueker,~G.S.: Two NP-complete problems in nonnegative integer programming. {\it Report No. 178}, Computer Science Laboratory, Princeton University, 1975.













\bibitem{oeis} {\it OEIS: The On-Line Encyclopedia of Integer Sequences.} \url{http://oeis.org/classic/index.html}. 

\bibitem{Pea05} Pearson,~D.: \lq\lq A polynomial-time algorithm for the change-making problem\rq\rq. {\it Operations Research Letters} \textbf{33}, 231--234 (2005).

\bibitem{PR21} P\'erez-Ros\'es,~H., M.~Bras-Amor\'os and J.M.~Serradilla-Merinero: \lq\lq Greedy routing in circulant networks\rq\rq. {\it Graphs and Combinatorics} \textbf{38}, 86 (2022). DOI: \url{https://doi.org/10.1007/s00373-022-02489-9}  

\bibitem{RoRo12} Robles-P\'erez,~A.M. and J.C.~Rosales: The Frobenius problem for numerical semigroups with embedding dimension equal to three. {\it Mathematics of Computation} 81 (2012), 1609--1617. DOI: \url{https://doi-org.sabidi.urv.cat/10.1090/S0025-5718-2011-02561-5}

\bibitem{RoGar09} Rosales,~J.C. and P.A.~Garc\'{\i}a-S\'anchez: {\it Numerical Semigroups}. Springer (2009).



\bibitem{Sha03} Shallit,~J.: \lq\lq What This Country Needs is an 18\cent \ Piece\rq\rq. {\it The Mathematical Intelligencer} \textbf{25}(2), 20--23 (2003).




\bibitem{Stra15} Strazzanti,~F.: \lq\lq Minimal genus of a multiple and Frobenius number of a quotient of a numerical semigroup\rq\rq. {\it International Journal of Algebra and Computation} \textbf{25}(6), 1043--1053 (2015).










\end{thebibliography}
\end{document}